\documentstyle[amscd,amssymb,verbatim]{amsart}
\newcommand{\disj}{\bigsqcup}
\newcommand{\Nm}{\operatorname{Nm}}

\newcommand{\GL}{\operatorname{GL}}

\newcommand{\Gal}{\operatorname{Gal}}

\newcommand{\SL}{\operatorname{SL}}

\newcommand{\lan}{\langle}
\newcommand{\ran}{\rangle}

\newcommand{\Tr}{\operatorname{Tr}}

\newcommand{\Th}{\Theta}

\newcommand{\eps}{\epsilon}

\numberwithin{equation}{section}

\newtheorem{prop}{Proposition}[section]

\newcommand{\Pf}{\noindent {\it Proof}}
\newcommand{\id}{\operatorname{id}}

\newcommand{\ov}{\overline}

\newcommand{\Aut}{\operatorname{Aut}}

\newcommand{\ra}{\rightarrow}

\newcommand{\la}{\lambda}

\newcommand{\Z}{{\Bbb Z}}
\newcommand{\Q}{{\Bbb Q}}
\newcommand{\F}{{\Bbb F}}

\newcommand{\La}{\Lambda}
\newcommand{\Ga}{\Gamma}

\newcommand{\sign}{\operatorname{sign}}

\newcommand{\ed}{\qed\vspace{3mm}}

\title{A new look at Hecke's indefinite theta series}
\author{A. Polishchuk}
\thanks{This work was partially supported by NSF grant}

\begin{document}

\maketitle

\bigskip

This note is devoted to the $q$-series of the form
$$\sum_{m\ge 0,n\ge 0}f(m,n)q^{Q(m,n)}-\sum_{m<0,n<0}f(m,n)q^{Q(m,n)}$$
where $Q$ is an indefinite quadratic form on $\Z^2$,
$f(m,n)$ is a doubly periodic function on $\Z^2$ such that
the sums of $f(m,n)q^{Q(m,n)}$ over all vertical and all horizontal
lines in $\Z^2$ vanish.
Some of these series appeared as coefficients in univalued triple Massey
products on elliptic curves
computed via homological mirror symmetry in \cite{P}. In particular,
in this context the condition of vanishing of sums over vertical and 
horizontal lines appears to be related to the standard necessary condition of 
the existence of triple Massey products (the vanishing of two double products).
In the present paper we generalize Theorem 3 of \cite{P} which relates
such series to the indefinite theta series
considered by Hecke in \cite{He1}, \cite{He2} 
(our approach is completely elementary and doesn't use the connection
with triple products on elliptic curves). 
The main consequence of this relation is the modularity of our
$q$-series. We also show that the problem of finding all linear relations
between our series 
is related to the study of orbits of actions of dihedral groups on
$(\Z/N\Z)^2$. 

\section{Main result}

\subsection{Hecke's indefinite theta series}
Let us recall the definition of these series. Let $K$ be a totally real
quadratic extension of $\Q$, i.e. $K$ is either a field of the form
$\Q(\sqrt{D})$ (where $D>0$)
or the algebra $\Q\oplus\Q$. We have the norm map
$\Nm:K\ra\Q$ (in case of $\Q\oplus\Q$ this is the product of components).
Let us denote by $C\subset K$ the set of elements with positive norm.
The cone $C$ is a union of two components and we define the function
$\sign:C\ra\pm 1$ which assigns value
$1$ (resp. $-1$) on totally positive (resp. negative) elements (in the case
of $\Q\oplus\Q$ ``total positivity" means positivity of both components).
Let us denote by $U_+(K)$ the subgroup of the multiplicative group
$K$ consisting of totally positive elements $k\in K^*$ with norm $1$
(in the case of $\Q\oplus\Q$ this is the group of elements $(r,r^{-1})$
where $r>0$).  Note that the group of $\Q$-linear automorphisms
of $K$ preserving $\Nm$ decomposes as follows:
$$\Aut_{\Q}(K,\Nm)={\pm \id}\times U_+(K)\times\Gal(K/\Q)$$
where $U_+(K)$ acts on $K$ by multiplication.
Let $\La\subset K$ be a lattice (i.e. a $\Z$-submodule of rank $2$),
$\La+c$ be a coset for this lattice (where $c\in K$).
Hecke's indefinite theta series is
$$\Th_{\La,c}=\sum_{\la\in(\La+c)\cap C/G}\sign(\la)q^{d\cdot \Nm(\la)}$$
where $G$ is the subgroup in $U_+(K)$ consisting of the elements
preserving $\La+c$, $d$ is a positive rational number such that
$d\Nm$ takes integer values on $\La+c$. Hecke proved that this series is
modular of weight $1$
for the subgroup $\Ga_0(n)\subset\SL_2(\Z)$ with some explicit level $n$.
\footnote{In the original definition of Hecke 
$\La$ was an ideal in the ring of integers, however, 
the same proof works for any lattice. Also, Hecke makes a 
concrete choice of $d$. For our
purposes it is more convinient to allow any $d$ such that $d\Nm$ takes
integer values on $\La+c$.}
Note that the elements of $U_+(K)$ preserving $\La$ are totally positive
units, hence, $G$ is an infinite
cyclic group. In particular, if we replace in the above definition $G$
by any infinite subgroup in $U_+(K)$ preserving $\La+c$ the resulting
series will be an integral multiple of $\Th_{\La,c}$.

\subsection{Formulation of the main theorem}

{\it
Let $Q(m,n)=am^2+2bmn+cn^2$ be a $\Q$-valued indefinite quadratic form
on $\Z^2$ (so $b^2>ac$) which is positive on the cone $mn\ge 0$ (i.e.
$a$, $b$ and $c$ are positive).
Let $f(m,n)$ be a doubly periodic complex-valued function on
$\Z^2$ (so $f(m+N,n)=f(m,n+N)=f(m,n)$ for some $N>0$).
Assume that for all $m_0$ and $n_0$ one has
$$\sum_{m\in\Z} f(m,n_0)q^{Q(m,n_0)}=
\sum_{n\in\Z}f(m_0,n)q^{Q(m_0,n)}=0$$
(i.e. all sums along horizontal and vertical lines are zero).
Assume also that $Q$ takes integer values on the support of $f$.
Then the series
$$\Th_{Q,f}=\sum_{m\ge 0,n\ge 0}f(m,n)q^{Q(m,n)}-
\sum_{m<0,n<0}f(m,n)q^{Q(m,n)}$$
is modular of weight $1$.

Moreover, the space of modular forms
of weight $1$ spanned by these series coincides with the space
generated by Hecke's indefinite theta series.
}

\subsection{Proof}
Our first task is to unravel the condition that the sums along
horizontal and vertical lines are zero.
Let us extend the function $f(m,n)$ from $\Z^2$ to $\Q^2$ by zero.
Then we claim that this condition is equivalent to the following
two identities:
$$f(m,n)=-f(-\frac{2b}{a}n-m,n),$$
$$f(m,n)=-f(m,-\frac{2b}{c}m-n).$$
Indeed, this follows from the fact that $Q$ restricted to
a vertical or horizontal line assumes each value at exactly two points
(sometimes coinciding, in which case the coefficient should be zero),
so in order for the sum to be zero the corresponding coefficients should
cancel out.
Let us consider the following two operators preserving $Q$:
$$A=\left(\matrix -1 & p \\ 0 & 1\endmatrix\right),$$
$$B=\left(\matrix 1 & 0 \\ r & -1\endmatrix\right).$$
where $p=-\frac{2b}{a}$, $r=-\frac{2b}{c}$.
Then the conditions on $f$ can be rewritten as
\begin{equation}\label{condition}
f(Ax)=f(Bx)=-f(x)
\end{equation}
for every $x\in\Q^2$. Let $S\subset\Z^2\subset\Q^2$ be the support of $f$.
We can assume that $f\neq 0$ so that $S$ is non-empty.
Let $\La=\{x\in\Q^2:\ S+x=S\}$. Since $f$ is doubly periodic, $\La$ is
a sublattice of $\Z^2$. On the other hand, both operators $A$ and $B$
preserve $S$, hence, they preserve $\La$. It follows
that $\Tr(AB)=-2+rp$ is an integer, i.e. $rp=\frac{4b^2}{ac}$ is an
integer.

Making the change of variables of the form
$m=m'/m_0$, $n=n'/n_0$, where $m_0$ and $n_0$ are positive integers
such that $\frac{m_0}{n_0}=\frac{a}{2b}$, we can always assume that
$a=2b$. Then the above condition will imply that both matrices
$A$ and $B$ have integer coefficients. In particular, we can consider
them acting on $(\Z/N\Z)^2$, where $N$ is the (double) period of $f$.
Let us denote by $G_N$ the subgroup of $\GL_2(\Z/N\Z)$ generated by
these two operators (by abuse of notation we will denote the corresponding
elements of $G_N$ also by $A$ and $B$).
Note that $A^2=B^2=1$, so $G_N$ is actually a dihedral group.
Now clearly the space of functions $f$ on $(\Z/N\Z)^2$ satisfying
the condition (\ref{condition}) is spanned by functions supported on
orbits of $G_N$ (and satisfying (\ref{condition})). Let
$O\subset (\Z/N\Z)^2$ be an orbit of $G_N$, $f$ be a function on
$O$ satisfying (\ref{condition}). In order for $f$ to be non-zero
the orbit $O$ should satisfy the following condition:
for every $x\in O$ one has $Ax\neq x$, $Bx\neq x$.
Let us call such orbit {\it admissible}. Conversely, it is easy to
see that for every admissible orbit $O$ there is a unique (up to a constant)
function $f$ on $O$ satisfying (\ref{condition}). Indeed, let
$\chi:G_N\ra\{\pm 1\}$ be the character defined by $\chi(A)=\chi(B)=-1$. 
Then the orbit is admissible if and only if $\chi$ is trivial on
the stabilizer subgroup of a point in $O$ (since every element $g\in G_N$
with $\chi(g)=-1$ is conjugate either to $A$ or to $B$).
Thus, for every admissible orbit $O=Gx$ we can define the function $f_O$ on
$O$ by setting $f_O(gx)=\chi(g)$ (up to a sign $f_O$ doesn't depend on $x$).
It suffices to deal with the series associated with such functions.
So, in the rest of the proof we will assume that $f$ is a doubly
periodic function on $\Z^2$ with values $\pm 1$ satisfying 
(\ref{condition}). Let $S\subset\Z^2\subset\Q^2$ be
the support of $f$. Then $S=S_1\cup S_{-1}$ where $S_1=f^{-1}(1)$,
$S_{-1}=f^{-1}(-1)$. Furthermore, we have $AS_1=BS_1=S_{-1}$.
Let $K$ be the quadratic extension of $\Q$ associated with the form $Q$.
If $D=b^2-ac$ is not a complete square then
$K$ is a real quadratic field $\Q(\sqrt{D})$, otherwise,
$K=\Q\oplus\Q$. 
The usual notation $x+y\sqrt{D}$ for elements of a real quadratic field
$K$ can be extended to the case when $D$ is a complete square and
$K=\Q\oplus\Q$. Namely, in this case we set
$x+y\sqrt{D}:=(x+y\sqrt{D},x-y\sqrt{D})$. We have
$$Q(m,n)=\frac{1}{c}[(bm+nc)^2-Dm^2]=\frac{1}{c}\Nm(bm+nc+m\sqrt{D}).$$
Thus, it makes sense to consider $\Z^2$ as a lattice in
$K$ via the map $(m,n)\mapsto (bm+nc+m\sqrt{D})$.
For two non-zero elements $k_1,k_2\in K$ let us denote 
$\lan k_1,k_2\ran=\Q_{>0}k_1+\Q_{>0}k_2$, $[k_1,k_2]=
\Q_{\ge 0}k_1+\Q_{\ge 0}k_2$, $\lan k_1,k_2]=\Q_{\ge 0}k_1+\Q_{>0}k_2$.
Using this notation we can write
\begin{align*}
\Th_{Q,f}=
&\sum_{\la\in S_1\cap [1,b+\sqrt{D}]}q^{\frac{\Nm(\la)}{c}}-
\sum_{\la\in S_1\cap \lan -1,-b-\sqrt{D}\ran}q^{\frac{\Nm(\la)}{c}}
-\sum_{\la\in S_{-1}\cap [1,b+\sqrt{D}]}q^{\frac{\Nm(\la)}{c}}+\\
&\sum_{\la\in S_{-1}\cap \lan -1,-b-\sqrt{D}\ran}q^{\frac{\Nm(\la)}{c}}.
\end{align*}
Let us extend the operators $A$ and $B$ from our lattice to $K$ by
$\Q$-linearity.
We have $B(1)=-1$, $B(b+\sqrt{D})=-b+\sqrt{D}$. Therefore, making the change
of variables $\la\mapsto B\la$ in the last two sums we get
$$\sum_{\la\in S_{-1}\cap [1,b+\sqrt{D}]}q^{\frac{\Nm(\la)}{c}}=
\sum_{\la\in S_1\cap [-1,-b+\sqrt{D}]}q^{\frac{\Nm(\la)}{c}},$$
$$\sum_{\la\in S_{-1}\cap \lan -1,-b-\sqrt{D}\ran}q^{\frac{\Nm(\la)}{c}}=
\sum_{\la\in S_1\cap \lan 1,b-\sqrt{D}\ran}q^{\frac{\Nm(\la)}{c}}.$$
Hence, we can rewrite $\Th_{Q,f}$ as follows:
$$\Th_{Q,f}=
\sum_{\la\in S_1\cap \lan b-\sqrt{D},b+\sqrt{D}]}q^{\frac{\Nm(\la)}{c}}-
\sum_{\la\in S_1\cap [-b+\sqrt{D},-b-\sqrt{D}\ran}q^{\frac{\Nm(\la)}{c}}.$$
Now it is easy to check that the operator
$AB:K\ra K$ coincides with multiplication by the element
$\frac{b+\sqrt{D}}{b-\sqrt{D}}$ of norm $1$. Therefore, we have
$$\Th_{Q,f}=\sum_{\la\in S_1\cap C/G}\sign(\la)q^{\frac{\Nm(\la)}{c}},$$
where $G$ is the infinite cyclic group generated by $AB$.
Note that the set $S_1$ is a union of a finite
number of cosets $(\La_1+x_i, i=1,\ldots,s)$ for the lattice
$\La_1=\{x\in K: S_1+x=S_1\}$. Furthermore, since $\La_1$ is preserved
by the action of $G$, there is a subgroup of finite index
$G_0\subset G$ preserving each of these cosets. Then we have
$$[G:G_0]\Th_{Q,f}=\sum_{\la\in S_1\cap C/G_0}
\sign(\la)q^{\frac{\Nm(\la)}{c}}=\sum_{i=1}^s
\sum_{\la\in(\La_1+x_i)\cap C/G_0}\sign(\la)q^{\frac{\Nm(\la)}{c}}.$$
Now each of the terms is a scalar multiple of Hecke's series.

Conversely, assume that we are
given a lattice $\La\subset K$ in a totally real
quadratic extension of $\Q$ and a coset $\La+c$.
Let $G\subset U_+(K)$ be the subgroup preserving $\La+c$.
Recall that $G$ is an infinite cyclic group.
Let $\eps$ be a generator of $G$.
Let us define the $\Q$-linear operators $A$ and $B$ on $K$ as follows:
$B(x)=-\ov{x}$ where $\ov{x}$ is the conjugate element to $x$ (in the
case $K=\Q\oplus\Q$ and $x=(x_1,x_2)$ one has $\ov{x}=(x_2,x_1)$),
$A(x)=-\eps\cdot\ov{x}$.
Note that $A^2=B^2=1$ while $\det A=\det B=-1$.
Let $k\in K$ be an eigenvector for $A$ with eigenvalue $-1$, so that
$\eps\ov{k}=k$. Changing $k$ by $-k$ if necessary we can assume that
$k$ is totally positive.
Then we have
$$\Th_{\La,c}=\sum_{\la\in(\La+c)\cap C/G}\sign(\la)q^{d\cdot \Nm(\la)}=\\
\sum_{\la\in (\La+c)\cap [k,\ov{k}\ran}  q^{d\cdot \Nm(\la)}-
\sum_{\la\in(\La+c)\cap\lan -k, -\ov{k}]} q^{d\cdot \Nm(\la)}.
$$
Note that we have $1\in\lan k,\ov{k}\ran$ since $k$ is totally positive.
Therefore, we can split each of the above sums into two
according to decompositions
$[k,\ov{k}\ran=[k,1]\disj \lan 1,\ov{k}\ran$,
$\lan -k,-\ov{k}]=\lan -k,-1\ran\disj [-1,-\ov{k}]$.
Making the change of variable $\la\mapsto B(\la)$ in the
sums over $\lan 1,\ov{k}\ran$ and over $[-1,-\ov{k}]$ we can rewrite
the above sum as follows:
$$\Th_{\La,c}=\sum_{\la\in S\cap ([1,k]\cup\lan -k,-1\ran)}f(\la)\sign(\la)
q^{d\cdot \Nm(\la)},$$
where $S=(\La+c)\cup B(\La+c)$, the function $f$ supported on $S$
is defined by
$$f(x)=\delta_{\La+c}(x)-\delta_{B(\La+c)}(x)$$
where $\delta_{I}$ is the characteristic function of the set
$I$. Note that since the operator $AB$ preserves $\La+c$ and
$(AB)B=B(AB)^{-1}$, it also preserves $B(\La+c)$, hence,
$f(ABx)=f(x)$. On the other hand, by definition $f(Bx)=-f(x)$.
Therefore, we also have $f(Ax)=-f(x)$. Now taking the coordinates
with respect to the basis $(1,k)$ as variables of summation
we see that the above series assumes the form
$$\sum_{(m,n)\in S,m\ge 0, n\ge 0}f(m,n)q^{Q(m,n)}-
\sum_{(m,n)\in S,m<0,n<0}f(m,n)q^{Q(m,n)}$$
where $S$ is a finite union of cosets with respect to some $\Z$-lattice in
$\Q^2$, $f$ is a periodic function on $S$ with the property that
sums of $f(m,n)q^{Q(m,n)}$ over all vertical and horizontal lines are zero.
It remains to change variables $(m,n)$ to $(Mm,Mn)$ where $MS\subset\Z^2$
to rewrite this series in the form we require.
\ed

\section{Remarks and examples}

\subsection{Linear relations}
The series $\Th_{Q,f}$ is often equal to zero.
It is an important open problem to formulate the necessary and
sufficient conditions for it to be zero. In other words,
the problem is to describe
all linear relations between such series for some basis in the
space of functions $f$ satisfying the assumptions of the main theorem.
We restrict ourself to
several observations. As above we assume that $p=-2b/a$ and $r=-2b/c$ are
integers, so that we have an action of operators $A$ and $B$ on
$\Z^2$ preserving the form $Q$. 
In the course of proof of the main theorem
we introduced the subgroup $G_N\subset\GL_2(\Z/N\Z)$ generated by these
two operators modulo $N$.
As we have seen above the space of functions on
$(\Z/N\Z)^2$ satisfying the condition (\ref{condition}) (further called
{\it admissible} functions) has a basis
$(f_O)$ enumerated by admissible $G_N$-orbits.
The change of variables $(m,n)\mapsto (-m,-n)$ shows that
$$\Th_{Q,f}=-\Th_{Q,f\circ[-1]},$$
where $f\circ[-1](m,n)=f(-m,-n)$.
Let us call an admissible orbit $O$ {\it symmetric} if $-O=O$, and
{\it asymmetric}
otherwise. Note that for an asymmetric orbit one has $O\cap -O=\emptyset$.
For every symmetric orbit $O$ the corresponding function $f_O$ is either
even or odd. We call a symmetric orbit $O$ {\it even} (resp. {\it odd})
if $f_O$ is even (resp. odd).
Now the above equation shows that for an even symmetric
orbit $O$ one has $\Th_{Q,f_O}=0$, while for an asymmetric orbit $O$
one has $\Th_{Q,f_O}=\pm\Th_{Q,f_{-O}}$ (the sign comes from the
sign ambiguity in the definition of $f_O$).

The action of the operator
$\tau:(m,n)\mapsto (n,m)$ gives some additional relations between
$\Th_{Q,f_O}$. Indeed, for any $Q$ we have
$$\Th_{Q,f}=\Th_{Q\circ \tau,f\circ \tau}.$$
If $Q\circ \tau=Q$ (i.e. $a=c$) then for every
admissible orbit $O$ we have $f_O\circ \tau=\pm f_{O'}$
for some other admissible orbit $O'$, hence
$\Th_{Q,f_O}=\pm\Th_{Q,f_{O'}}$.
In particular, if $f_O\circ \tau=-f_O$ then $\Th_{Q,f_O}=0$.

Finally, we can make the changes of variables $(m,n)\mapsto(t_1 m,t_2 n)$,
where $t_1$ and $t_2$ are positive rational numbers, in the case
when this transformation sends the support of $f$ into $\Z^2$.
This transformation will always change the form $Q$ (unless
$t_1=t_2=1$). However, combining it with the operator $\tau$ with respect
to the new variables we can derive more linear relations for fixed $Q$
(generalizing the above relations for the case $a=c$).
Namely, assume that $c/a=t^2$ for some positive rational number
$t$. Then the operator
$$\tau_{t}:(m,n)\mapsto (t n,t^{-1}m)$$
preserves $Q$ and satisfies $\tau_{t}^2=1$, $\tau_{t}A=B\tau_{t}$.
In particular if $f$ is an admissible function such that
$\tau_{t}$ sends the support of $f$ into $\Z^2$ then
$f\circ \tau_{t}$ is also admissible (perhaps with a different double
period) and we have
$\Th_{Q,f\circ \tau_{t}}=\Th_{Q,f}$.

We were not able to find any other linear relations between
the series $(\Th_{Q,f})$ for fixed $Q$.
However, at present we are far from proving
that these are all relations. Even the non-vanishing of $\Th_{Q,f_O}$
for odd symmetric and for asymmetric admissible $G_N$-orbits
(in the case when $a/c$ is not a square in $\Q$) is still an open problem.
Note that some non-vanishing results were proven in \cite{P} using
homological mirror symmetry.

In Hecke's paper one can find the following vanishing
condition for an indefinite theta series $\Th_{\La,c}$ (see \cite{He2},
Satz 1): if there exists a totally negative element $\delta\in K^*$
with $\Nm(\delta)=1$ such that
$\delta(\La+c)=\La+c$ then $\Th_{\La,c}=0$. Let us show
that this vanishing is actually one of the linear relations considered
above. We will use the notation introduced in the proof of the main theorem.
Let $(Q,f)$ be the data constructed in the second half of the proof
so that $\Th_{\La,c}=\Th_{Q,f}$.
First of all, notice that $\delta^2\in G$, hence $\delta^2=\eps^n$
for some integer $n$. Changing $\delta$ by a power of $\eps$ we can
assume that either $\delta^2=1$ or $\delta^2=\eps$. In the former
case $\delta=-1$ so one has $f\circ[-1]=f$. In the latter case
we have $\eps\ov{\delta}=\delta$, so rescaling $k$ we can assume
that $k=-\delta$. It is easy to see that the operator
$\delta B:x\mapsto -\delta \ov{x}$
preserves $\Nm$ and switches $1$ and $k$, as well as $\La+c$ and
$B(\La+c)$.
Hence, the transposition $\tau:(m,n)\mapsto (n,m)$ preserves $Q$
and satisfies $f\circ \tau=-f$. A different choice of $k$ would lead
to a similar relation with $S$ replaced by $\tau_{t}$.

\subsection{Symmetric orbits}

Henceforward, operators $A$ and $B$ are always considered modulo $N$.
In the situation when the subgroup $G_N\subset\GL_2(\Z/N\Z)$ ($N>2$)
contains the matrix $-\id$ every orbit is symmetric. Furthermore,
since the character $\chi:G_N\ra\{\pm 1\}$ defined by $\chi(A)=\chi(B)=-1$
coincides with $\det|_{G_N}$, we have $\chi(-\id)=1$, hence every
orbit is even. Thus, we get $\Th_{Q,f}=0$ for all admissible
$f$. The following proposition gives a criterion allowing to recognize
this situation in the case when $N$ is an odd prime.

\begin{prop}\label{opp}
Assume that $N$ is an odd prime.
Then $-\id\in G_N$ if and only if $rp\mod N$ is of the
form $2+\la+\la^{-1}$ where $\la$ is an element of
even order in $\F_{N^2}^*$ ($\F_{N^2}$ is
the finite field of cardinality $N^2$). The number of such residues
modulo $N$ is equal to $N-\frac{n_1+n_2}{2}$ where
$n_1$ (resp. $n_2$) is the maximal odd divisor of $N-1$ (resp. $N+1$).
\end{prop}

\Pf . Since $G_N\cap\SL_2(\Z/N\Z)$ is generated by $AB$ the condition
$-\id\in G_N$ is equivalent to $(AB)^n=-\id$ for some $n$. We have
$\Tr(AB)=rp-2$, $\det(AB)=1$, so the eigenvalues $\la_1,\la_2$
of $AB$ are roots of the equation
$$\la^2-(rp-2)\la+1=0.$$
Assume first that $\la_1=\la_2$. Then either $rp=4$ or $rp=0$.
In the former case $\la_1=\la_2=
1$, hence, no power of $AB$ equals $-\id$.
In the latter case one can easily check that $(AB)^N=-\id$.
On the other hand, $0$ can be represented in the form $2+\la+\la^{-1}$
for $\la=-1$.

Now assume that $\la_1\neq\la_2$. Then the condition $(AB)^n=-\id$
is equivalent to $\la_1^n=-1$, i.e. $\la_1$ has even order in
the multiplicative group of $\ov{\F}_N$. It remains to notice that
$\la_1\in\F_{N^2}^*$ and that we have $rp-2=\la_1+\la_1^{-1}$.

To compute the number of such residues modulo $N$ we note that the condition
$\la+\la^{-1}\in\F_N$ means that either $\la^{N-1}=1$ or $\la^{N+1}=1$.
The number of elements $\la$ of even order such that $\la^m=1$ (where $m$ is
either $N-1$ or $N+1$)
is equal to $m-n$ where $n$ is the maximal odd divisor of $m$.
Therefore, the number of elements in $\F_N$ of the form $\la+\la^{-1}$
is equal to 
$$1+\frac{(N-1)-n_1-1}{2}+\frac{(N+1)-n_2-1}{2}=N-\frac{n_1+n_2}{2}.$$
\ed

Taking in the above proposition $\la$ to be $-1$, $\zeta_4$ and $\zeta_6$
(where $\zeta_l$ is a primitive root of unity of order $l$)
we get $rp\equiv 0\mod(N)$, $rp\equiv 2\mod(N)$ and 
$rp\equiv 3\mod(N)$ respectively.
On the other hand, we claim that if $rp\equiv 1\mod(N)$ or
$rp\equiv 4\mod(N)$ then $-\id\not\in G_N$.
Indeed, the equation $4=2+\la+\la^{-1}$ has the only solution $\la=1$
while the solutions of the equation $1=2+\la+\la^{-1}$ are roots of unity
of order $3$. These are the only cases of the above criterion which are
independent of $N$.
Here are the lists of values of $rp\mod(N)$
such that $-\id\not\in G_N$ for small odd primes $N$:

\noindent
$N=3$: $rp\equiv 1\mod(3)$;

\noindent
$N=5$: $rp\equiv 1,4\mod(5)$;

\noindent
$N=7$: $rp\equiv 1,4\mod(7)$;

\noindent
$N=11$: $rp\equiv 1,4,5,9\mod(11)$.

\noindent
$N=13$: $rp\equiv 1,4,9,10,12\mod(13)$.

Our last general observation is that in the case when $N$ is an odd prime,
all symmetric $G_N$-orbits have the same parity, i.e. 
they are either all odd or all even.

\begin{prop}\label{symodd}
Assume that $N$ is an odd prime. Then either $-\id\in G_N$
or every symmetric $G_N$-orbit is odd.
\end{prop}

\Pf . Assume that there exists a non-zero vector
$v\in (\Z/N\Z)^2$ and an element
$g\in G_N$ such that $gv=-v$ and $\det(g)=1$.
Then both eigenvalues of $g$ are $-1$, hence, $g^N=-\id$. 
\ed

\subsection{Examples}
In all examples below we assume that 
$a$, $c$, $p=-\frac{2b}{a}$ and $r=-\frac{2b}{c}$ are integers
($b$ is a half-integer).
Note that we are interested only in the cases when $G_N$
doesn't contain $-\id$. In particular, if $N$ is an odd prime we can
assume that $rp\not\equiv 0\mod(N)$. In this case 
the conjugacy class of the subgroup $G_N\subset\GL_2(\Z/N)^2$
depends only on $rp\mod(N)$. For instance, if $rp\equiv 1\mod(N)$ then
$G_N$ is isomorphic to the permutation group $S_3$. 
In examples 1 and 2 below we consider in details
cases $N=3$ and $N=5$. It turns out that in these cases
all admissible orbits are symmetric (they are automatically odd by
proposition \ref{symodd}). The simplest example
of an asymmetric admissible orbit (for prime $N$) occurs for $N=7$
(see example 3 below).

\noindent
1. $N=3$, $rp\equiv 1\mod(3)$.
Then there is a unique admissible orbit: the orbit of $(1,0)$.
For $r\equiv p\equiv 1\mod(3)$
(resp. $r\equiv p\equiv -1\mod(3)$)
the corresponding admissible function is
$f(m,n)=\chi_3(m+n)$ (resp. $f(m,n)=\chi_3(m-n)$)
where $\chi_3$ is the non-trivial
Dirichlet character modulo $3$ such that $\chi_3(\pm 1)=\pm 1$.
Let us assume that $a\le c$
(we can always achieve this using the transformation
$(m,n)\mapsto (n,m)$ if necessary).
Then we have
$$\Th_{Q,f}\equiv q^a+\chi_3(r)q^c\mod(q^{a+1}).$$
It follows that this theta series doesn't vanish
unless $r\equiv -1\mod(3)$ and $a=c$. In the latter case
we have $Q(n,m)=Q(m,n)$ while $f(n,m)=-f(m,n)$ so that $\Th_{Q,f}=0$.

\noindent
2. $N=5$.

\noindent
(a) $rp\equiv 1\mod(5)$. In this case there
are two distinct admissible orbits: the orbit of $(1,0)$ and
the orbit of $(2,0)$. It is easy to see that unless $a=c$
the corresponding two theta functions $\Th_{Q,f_1}$ and $\Th_{Q,f_2}$
are linearly independent. More precisely, the initial terms of these
series look as follows (in (i) and (ii) we assume that $a\le c$):

(i)$p\equiv r\equiv 1(5)$:
$$\Th_{Q,f_1}\equiv q^a+q^c\mod(q^{a+1}),\
\Th_{Q,f_2}\equiv q^{4a}+q^{4c}\mod(q^{4a+1}).$$

(ii)$p\equiv r\equiv -1(5)$:
$$\Th_{Q,f_1}\equiv q^a-q^c\mod(q^{a+1}),\
\Th_{Q,f_2}\equiv q^{4a}-q^{4c}\mod(q^{4a+1}).$$

(iii)$p\equiv 2(5)$, $r\equiv -2(5)$:
$$\Th_{Q,f_1}\equiv q^a-q^{4c}\mod(q^{\min(a,4c)+1}),\  
\Th_{Q,f_2}\equiv q^{c}-q^{4a}\mod(q^{\min(4a,c)+1}).$$
Furthermore, in the case $a=4c$ we have
$$\Th_{Q,f_1}\equiv q^{9c}\mod(q^{9c+1})$$ while in the case
$c=4a$ we have
$$\Th_{Q,f_2}\equiv q^{9a}\mod(q^{9a+1}).$$

If $a=c$ then in the case (ii) we have $\Th_{Q,f_1}=\Th_{Q,f_2}=0$
while in the case (iii) we have $\Th_{Q,f_2}=\Th_{Q,f_1}$.

\noindent
(b) $rp\equiv -1\mod(5)$. In this case $AB$ has order $5$
but there are still two admissible orbits:
the orbit of $(1,0)$ and the orbit of $(2,0)$.\footnote{In this case
$A$ and $B$ have a common invariant vector which allows to have
bigger admissible orbits than in case (a).}
The analysis of the initial terms of these series
(very similar to the case (a)) implies that the corresponding
two theta series are linearly independent unless $a=c$.

\noindent
3. $N=7$, $r\equiv p\equiv 1\mod(7)$. There are $5$ admissible orbits:
$3$ symmetric orbits and $2$ asymmetric orbits.
The symmetric orbits are $O_1=G_N\cdot (1,0)$, $2\cdot O_1$,
and $3\cdot O_1$. The asymmetric orbits are $O_2=G_N\cdot (1,3)$ and
$-O_2$. Using the relation $\Th_{Q,f_{-O_2}}=-\Th_{Q,f_{O_2}}$ we can
exclude the orbit $-O_2$ from our consideration. The initial
terms of the remaining $4$ theta series look as follows
(assuming that $a\le c$)
$$\Th_{Q,f_{O_1}}\equiv q^a+q^c\mod (q^{a+1}),$$
$$\Th_{Q,f_{2O_1}}\equiv q^{4a}+q^{4c}\mod (q^{4a+1}),$$
$$\Th_{Q,f_{3O_1}}\equiv q^{9a}+q^{9c}\mod (q^{9a+1}),$$
$$\Th_{Q,f_{O_2}}\equiv q^{9a+c+6b}+q^{a+9c+6b}\mod (q^{9a+c+6b+1}).$$
This immediately implies that they are linearly independent.

\noindent
4. The indefinite theta series considered in Theorem 2 of \cite{P}
correspond to the following situation. Let us assume that
$\frac{b}{a}$ and $\frac{b}{c}$ are integers (not just half-integers
as before). In this case the discriminant $D=b^2-ac$ is divisible by
$ac$. We are going to take $N=\frac{4D}{ac}$. 
Let $s_1$ and $s_2$ be arbitrary odd numbers.
It is easy to check 
that the $G_N$-orbit of the element 
$$v_{s_1,s_2}=(\frac{b}{a}s_2-s_1,\frac{b}{c}s_1-s_2)\in(\Z/N\Z)^2$$
is admissible and consists of four elements which are congruent
to $v_{s_1,s_2}$ modulo $N/2$. On the other hand, 
if $l$ divides $\frac{b}{a}+1$ and $\frac{b}{c}+1$ then 
$\frac{D}{ac}$ is divisible by $l$ and the $2l$-torsion element  
$v_l=\frac{2D}{lac}(1,1)\in(\Z/N\Z)^2$ is $G_N$-invariant.
The series considered in \cite{P} correspond to the orbits
of the elements $v_{s_1,s_2}+t\cdot v_l$  where $t\in \Z$
(these orbits depend only on $t\mod(l)$).

\end{document}